\documentclass{specgeom}
\usepackage{amscd,amssymb}

\newtheorem{theorem}{Theorem}[section]
\newtheorem{lemma}[theorem]{Lemma}
\newtheorem{corollary}[theorem]{Corollary}
\newtheorem{proposition}[theorem]{Proposition}

\theoremstyle{definition}
\newtheorem{definition}[theorem]{Definition}

\theoremstyle{remark}
\newtheorem{remark}[theorem]{Remark}

\numberwithin{equation}{section}

\def\CC{\mathbb{C}}
\def\RR{\mathbb{R}}

\def\lra{\longrightarrow}

\def\a{\alpha}
\def\ga{\gamma}
\def\la{\lambda}
\def\e{\varepsilon}
\def\si{\sigma}
\def\ph{\varphi}
\def\om{\omega}
\def\Ga{\Gamma}
\def\Om{\Omega}
\def\La{\Lambda}

\DeclareMathOperator\im{Ran} \DeclareMathOperator\Ker{Ker} \DeclareMathOperator\Coker{Coker}
\DeclareMathOperator\supp{supp} \DeclareMathOperator\const{const} \DeclareMathOperator\tr{tr}
\DeclareMathOperator\Tr{Tr} \DeclareMathOperator\ind{ind} \DeclareMathOperator\id{id}
\DeclareMathOperator\cch{ch}

\newcommand{\op}[1]{\operatorname{#1}}

\def\cA{{\mathcal{A}}}
\def\cD{{\mathcal{D}}}
\def\cH{{\mathcal{H}}}
\def\cK{{\mathcal{K}}}
\def\cL{{\mathcal{L}}}

\def\gW{{\wt{\cH}}}
\def\ov{\overline}
\def\wt{\widetilde}
\def\wh{\widehat}
\def\pa{\partial}

\newcommand\pd[2]{\frac{\pa #1}{\pa #2}}

\let\rom\textup

\newcommand{\BL}{\biggl}
\newcommand{\BR}{\biggr}
\newcommand{\bl}{\bigl}
\newcommand{\br}{\bigr}
\newcommand{\norm}[1]{\left\| #1 \right\|}

\newcommand{\dirac}{\setminus\!\!\!\!D}

\newcommand{\circM}{{M}^{\circ}}
\newcommand{\circU}{{U}^{\circ}}
\newcommand{\wA}{\wh{A}}
\newcommand{\wb}{\wh{b}}
\newcommand{\wP}{\wt{P}}
\newcommand{\wDo}{\bA\biggl(-i\frac\pa{\pa t}\biggr)}
\newcommand{\wAo}{\bA\biggl(-i\frac\pa{\pa t}\biggr)}
\newcommand{\Tau}{\mathrm{T}}
\def\bA{\mathbf{A}}
\def\bhA{\skew{-2}\wh{\mathbf{A}}}

\def\chhi#1{\chi_{{}_{\scriptstyle{#1}}}}

\newcommand{\dbb}{\bar\pa_b}

\begin{document}

\title[GUILLEMIN TRANSFORM AND TOEPLITZ REPRESENTATIONS]
{Guillemin Transform and Toeplitz Representations\\ for Operators on Singular Manifolds}

\author[V.~E.~NAZAIKINSKII]{V.~E.~Nazaikinskii}
\address{Institute for Problems in Mechanics, RAS, Moscow, Russia}
\email{nazaikinskii@yandex.ru}
\thanks{ The first, third and fourth authors thank the RFBR grants nos.~02-01-00118 and~03-02-16336, and The Swedish  Royal Academy of Sciences for financial support and Chalmers University of Technology for hospitality.}

\author[G.~ROZENBLUM]{G.~Rozenblum}
\address{Chalmers University of Technology , G\"oteborg, Sweden}
\email{grigori@math.chalmers.se}

\author[A.~Yu.~SAVIN]{A.~Yu.~Savin}
\address{Independent University of Moscow, Moscow, Russia}
\email{antonsavin@mail.ru}

\author[B.~Yu.~STERNIN]{B.~Yu.~Sternin}
\address{Independent University of Moscow, Moscow, Russia}
\email{sternin@mail.ru}

\subjclass[2000]{ Primary 46L80; Secondary 58L80, 46L87, 58J40}

\keywords{Elliptic theory on manifolds with singularities, index
formulas, pseudodifferential quantization, Toeplitz quantization,
Szeg\"o--Calder{\'o}n projection, Guillemin transform}

\begin{abstract}
A new approach to the construction of index formulas for elliptic
operators on singular manifolds is suggested on the basis of
$K$-theory of algebras and cyclic cohomology. The equivalence of
Toeplitz and pseudodifferential quantizations, well known in the
case of smooth closed manifolds, is extended to the case of
manifolds with conical singularities. We describe a general
construction that permits one, for a given Toeplitz quantization
of a $C^*$-algebra, to obtain a new equivalent Toeplitz
quantization provided that a resolution of the projection
determining the original quantization is given.
\end{abstract}

\maketitle

\section*{Introduction}

In  recent years, the problem of finding index formulas for
elliptic pseudodifferential operators on singular manifolds has
been studied in numerous papers (e.g., see
    Plamenevsky-Rozenblum~\cite{PlRo4}, Melrose-Nistor~\cite{MeNi1}, Rozenblum~\cite{Roz1},
     Schulze-Sternin-Shatalov~\cite{ScSS11}, Nazaikinskii-Sternin~\cite{NaSt3}). Despite considerable
progress in this direction, the situation is far from being clear
yet. Indeed, some of the formulas known in the literature fail to
express the index via the principal symbol treated as an element
of an appropriate Calkin algebra (e.g.,~\cite{MeNi1}), in other
formulas, separate terms are not homotopy invariant (e.g., see
Fedosov-Schulze-Tarkhanov~\cite{FST1}), and the few formulas that
combine both desirable properties (e.g., see~\cite{ScSS11,NaSt3})
are valid only for an important but rather narrow class of
operators satisfying certain symmetry conditions.

This situation is primarily caused by the complicated symbol
structure for pseudodifferential operators on singular manifolds.
These symbols consist of components (in general, operator-valued)
corresponding to the strata of the manifold and satisfying certain
matching conditions for the adjacent strata see, e.g.
Schulze~\cite{Schu1}. In a number of existing results, the index
of the operator is expressed as a sum of contributions from  these
symbol components. These contributions are usually noninteger and
lack homotopy invariance, which is not at all surprising in the
presence of the matching conditions. It is also difficult to
assign to them  any straightforward topological or algebraic
meaning. One can ensure their homotopy invariance only by severely
restricting the class of operators to be considered (say, by
imposing some symmetry conditions). A detailed analysis of
symmetry conditions and their role in obtaining invariant index
formulas can be found in Savin-Schulze-Sternin~\cite{SaScS1} and
in \cite{NaSt3}, and we do not dwell on the topic here.

In the present paper, we propose another approach to the construction of index formulas on
manifolds with singularities. This approach is based on $K$-theory of algebras and cyclic
cohomology. We take a slightly different viewpoint as to what a ``good'' index formula must be.
Instead of trying to use topological invariants of separate components of the symbol (which is
hopeless due to results of ~\cite{SaScS1}), we consider the symbol as a whole, that is, as an
element of an appropriate symbol algebra $\cA$. Moreover, instead of topological objects, one
naturally deals with algebraic objects like the Chern--Connes character viewed as an element of
the cyclic cohomology group of the symbol algebra. In the abstract framework, the scheme is well
known.

If the algebra $\cA$ is separable, then every quantization of
$\cA$ (in particular, the pseudodifferential quantization, which
is of interest to us) is equivalent\footnote{Quantizations are
said to be \textit{equivalent} if they determine the same element
in the group $K^1(\cA)$ (e.g.,~see Blackadar~\cite{Bla1}).} to
some (generalized) Toeplitz quantization
\begin{equation}\label{toep}
\begin{CD}
    \wt\tau(a)=Pa\,\,:\cH\hookrightarrow \wt\cH@>{a}>>
    \wt\cH@>P>>\cH
\end{CD}
\end{equation}
in the Hilbert space $\cH=\im(P)$ defined as the range of an orthogonal projection $P$ in some
Hilbert $\cA$-module $\wt\cH$; here $P$ is assumed to almost commute with the action of $\cA$.

Under the additional condition that the commutators of elements of
the algebra with the projection are not only compact but also
belong to an appropriate von Neumann--Schatten class, the Toeplitz
representation permits one to write out analytic index formulas
(see Connes~\cite{Con1}). Namely, if $\cA_\infty\subset\cA$ is a
dense local subalgebra of $\cA$ such that the restriction of the
Toeplitz quantization to $\cA_\infty$ is
$p$-summable,\footnote{That is, the commutator $[a,P]$ belongs to
the $p$th von Neumann--Schatten class $\mathfrak{S}_p$ for every
$a\in\cA_\infty$.} then for an arbitrary elliptic element
$\mathfrak{a}\in \op{M}(\cA_\infty)$ one has the index formula
\begin{equation}
\label{index}
 \ind(\tau(\mathfrak{a}))=\frac1{2^{N-1/2}\sqrt{i}}\Gamma\BL(\frac N2+1\BR)
 \times\bigl(\cch[\tau]\otimes\tr\bigr)
 \bigl(\underbrace{\mathfrak{a}^{-1},\mathfrak{a},\mathfrak{a}^{-1},
 \mathfrak{a},  \ldots,\mathfrak{a}^{-1},
 \mathfrak{a}}_{\text{$N+1$ arguments}}\bigr),
\end{equation}
where $N>p$ is odd, $\tr$ is the matrix trace, and the \textit{Chern--Connes character}
$\cch[\tau]$ of the quantization $\tau$ is given by the formula
\begin{multline}
\label{chern}
 \cch[\tau](a_0,a_1,\ldots,a_N)\\=\sqrt{2i}
 (-1)^{N(N-1)/2}\Gamma\BL(\frac N2+1\BR)^{-1}\Tr\{a_0[P,a_1][P,a_2]\dotsm
[P,a_N]\}.
\end{multline}
(Here $\Tr$ is the operator trace in $\wt\cH$.)

To apply this formalism to manifolds with singularities, we do the following. For the case of a
smooth compact closed manifold $M$, there is a well-known Toeplitz quantization
\begin{equation}\label{toep1}
    a(x,p)\longmapsto \Pi a:\im(\Pi)\lra \im(\Pi),\quad
    \im(\Pi)\subset L^2(S^*M)
\end{equation}
of the symbol algebra $\cA\equiv C(S^*M)$ with the help of the
Szeg\"o-Calder{\'o}n projection $\Pi$ in the space $L^2(S^*M)$ of
square integrable functions on the cosphere bundle $S^*M$, see
Guillemin, ~\cite{Gui1}, Boutet de Monvel~\cite{Bout3}, Boutet de
Monvel-Guillemin~\cite{BoGu1}. This construction was explained
in~\cite{Gui1} as showing that pseudodifferential operators are
Toeplitz operators in disguise. In the present paper, we
generalize this quantization and construct a Toeplitz quantization
equivalent to the standard pseudodifferential quantization for
manifolds with conical singularities. Thus pseudodifferential
operators on a manifold with conical singularities are Toeplitz
operators in disguise, too. (This Toeplitz quantization acts in a
more complicated space, which, however, turns into $L^2(S^*M)$ if
there are no singularities.) In the smooth case, the equivalence
(and even the almost isomorphism) of the Toeplitz quantization to
the pseudodifferential quantization\footnote{In~\eqref{PSD}, the
continuation of the symbol on $T^*M$ as a homogeneous function of
degree zero in the fibers is assumed.}
\begin{equation}\label{PSD}
 a(x,p)\longmapsto \wh a=a\BL(x,-i\pd{}x\BR)
 :L^2(M)\lra L^2(M)
\end{equation}
is given by the so-called Guillemin transform
\begin{equation}\label{gui}
    \Tau:L^2(M)\lra L^2(S^*M)
\end{equation}
(see the cited papers), and we construct an analog of this
transform for manifolds with singularities, thus proving the
equivalence (and even an almost isomorphism) of quantizations in
this case. This naturally results in an index formula of the
form~\eqref{index} for elliptic operators on manifolds with
conical singularities. This algebraic index formula has several
advantages:
\begin{itemize}
  \item first, it expresses
the index via the principal symbol alone and is homotopy invariant;
  \item second, it is valid for arbitrary elliptic symbols,
  not just for symbols satisfying some symmetry conditions;
  \item third, most importantly, it is expressed in terms of the cyclic
  cohomology class $\cch[\tau]$ of the algebra $\cA_\infty$. It is
  the cyclic cohomology of the symbol algebra that replaces the
  homology of a manifold as one passes from the algebra
  of functions on a smooth manifold to more general symbol algebras.
\end{itemize}

\medskip

A serious disadvantage of this formula (and in general of the
index formula~\eqref{index} for Toeplitz quantizations given by
projections of general form) is that it fails to be local. At the
same time, if the chosen projection is the positive spectral
projection of an unbounded local operator, then a different
representation of the index cocycle~\eqref{chern} is possible,
which results in more traditional, local index formulas. (See
Connes-Moscovici~\cite{CoMo1}.) In this connection, the construction
of a local index formula involves the natural problem of finding
an unbounded local operator such that the Toeplitz quantization
generated by its positive spectral projection is equivalent to the
pseudodifferential quantization on a manifold with conical
singularities.

In the smooth case, the desired operator is the Dirac operator. More precisely, the natural
self-adjoint Dirac operator $\dirac$ acting on sections of the spinor bundle on the
odd-dimensional contact manifold $S^*M$ is associated with the almost complex structure on the
distribution of contact hyperplanes.  Then the index of a pseudodifferential operator
\begin{equation*}
 \wh a:C^\infty(M,\CC^N)\lra C^\infty(M,\CC^N)
\end{equation*}
proves to be equal to the index of the Toeplitz operator
constructed from the symbol~$a$ and the positive spectral
projection of the Dirac operator~(see Baum-Douglas~\cite{BaDo1}):
\begin{equation}\label{bado}
    \ind\wh a=\ind\bl\{P_+(\dirac)\otimes a:
    \im(P_+(\dirac))\otimes\CC^N
    \lra\im(P_+(\dirac))\otimes\CC^N\br\},
\end{equation}
where $\pi:S^*M\lra M$ is the natural projection. This equality of
indices shows that (at least, modulo torsion elements in the group
$K^1(\cA)$) the quantization given by the positive spectral
projection $P_+(\dirac)$ of the Dirac operator is equivalent to
the pseudodifferential quantization and hence to the quantization
with the help of the Szeg\"o--Calder{\'o}n projection. There are
$K$-theoretic proofs of this fact. (One of then is based on the
Atiyah--Singer theorem, e.g., see Kaminker~\cite{Kam1,Kam2}, and
another, by Baum-Douglas-Taylor~\cite{BDT1}, uses some facts
concerning the $\bar\pa$-Neumann problem.)

Both these approaches to the proof of~\eqref{bado} encounter
serious difficulties in the case of manifolds with singularities.
The main difficulty is that so far one does not have even a
hypothetical candidate for the ``Dirac operator'' (i.e., an
operator for which~\eqref{bado} holds) for the case of manifolds
with conical singularities.

A prototype for a possible generalization of formula~\eqref{bado}
to the conical case is the proof of~\eqref{bado} (as a special
case of more general index formulas for Toeplitz operators and
quantized contact transformations on general contact manifolds)
due to Epstein, Melrose, and Mendoza (see~\cite{EpMe1} and
also~\cite{HHH}). Their proof \textit{relates the
Szeg\"o--Calder{\'o}n projection and the positive spectral
projection of the Dirac operator by a finite chain of
transformations} and hence is of interest to us as a model for the
possible definition of the Dirac operator in the conical case.
Indeed, we have already constructed a counterpart of the
Szeg\"o--Calder{\'o}n projection for manifolds with conical
singularities and proved the equivalence of the Toeplitz
quantization and the pseudodifferential quantization. The
Epstein--Melrose--Mendoza construction comprises three steps (we
give more detail in the Appendix):
\begin{enumerate}
    \item constructing a resolution of the Szeg\"o--Calder{\'o}n projection (in the case of the
          cosphere bundle of a smooth manifold, the resolution is given by the
          Kohn--Rossi complex, e.g., see Folland-Kohn~\cite{FoKo1});
    \item proving that the Toeplitz quantization associated with the positive
          spectral projection of the operator
          $G$ obtained as the standard roll-up of this resolution produces the operator
          with the same index as the Toeplitz quantization associated with the
          Szeg\"o--Calder{\'o}n
          projection;
    \item deforming the self-adjoint operator $G$ to the Dirac operator $\dirac$.
\end{enumerate}

In the present paper, we show that step~2 of this construction is a special case of a general
assertion that permits one, starting from a Toeplitz quantization of an arbitrary $C^*$-algebra
$\cA$, to obtain new equivalent (i.e., defining the same class in $K^1(\cA)$) Toeplitz
quantizations associated with positive spectral projections of self-adjoint operators.

The construction goes as follows: for a projection $P$ almost
commuting  with the action of $\cA$, one takes a finite-length
resolution in Hilbert $\cA$-modules. Rolling up the resolution in
the standard way and adding the projection itself and the
projection onto the cokernel in the last term of the resolution,
one obtains a self-adjoint operator in the direct sum of all
Hilbert modules occurring in the resolution. The positive spectral
projection of this self-adjoint operator is the desired new
projection.

In conclusion, we note that the so far open problem of defining the Dirac operator $\dirac$ in
the conical case is undoubtedly of interest in that such an operator is a ``fundamental cycle''
of the algebra $\cA$ (the noncommutative manifold) and, in particular, defines the Poincar\'e
duality on it:
\begin{align*}
    K_1(\cA)&\overset{\simeq}\lra K^0(\cA),\\
    a&\longmapsto P_+(\dirac)a:\im P_+(\dirac)\lra\im P_+(\dirac).
\end{align*}
In the smooth case ($\cA=C(S^*M)$), this is reduced to the isomorphism
\begin{equation*}
    K^1(S^*M)\simeq K_0(S^*M).
\end{equation*}

\medskip

The paper is organized as follows. In the first
section, we construct a Toeplitz representation of the
pseudodifferential quantization on manifolds with
conical singularities. Section~2 contains the
above-mentioned construction of equivalent
quantizations from resolutions of projections. First
we consider the case in which the resolution is formed
by bounded operators and then the case in which the
resolution is formed by unbounded operators. In the
Appendix, we briefly describe the construction
of~\cite{EpMe1,HHH} for the special case in which the
contact manifold in question is the cosphere bundle of
a smooth compact manifold without boundary.

Throughout the paper, we assume that the reader is acquainted with the notions and definitions of
$K$-theory of operator algebras and related topics (e.g., see~\cite{Bla1,Con1}). On the other
hand, the necessary definitions and facts of the theory of (pseudo)differential operators on
manifolds with conical singularities are given in Subsection~1.1 together with appropriate
bibliographic references.

\section{The Toeplitz representation of pseudodifferential operators on a manifold with conical
singularities}

\subsection{Manifolds with conical singularities and $\Psi$DO}

We consider pseudodifferential operators on a compact manifold $M$
with isolated conical singularities. We always assume that there
is only one conical point $\a$; extending the results to the
general case encounters no difficulties. In this subsection we
briefly list main definitions and facts.  Details can be found,
e.g., in Schulze~\cite{Schu1}, Egorov-Schulze~\cite{EgSc1},
Schulze-Sternin-Shatalov~\cite{ScSS20} and, as it concerns the
cylindrical representation for cone-degenerate pseudodifferential
operators, also in Nazaikinskii-Sternin~\cite{NaSt12}.
\begin{definition}
A \textit{manifold $M$ with conical singularity} $\a$ is a compact Hausdorff topological space
$M$ with a distinguished point $\a$ such that the set $\circM=M\setminus\{\a\}$ is equipped with
the structure of a $C^\infty$ manifold compatible with the topology and there is a given
homeomorphism
\begin{equation}\label{zwezda}
\ph:U\simeq K_{\Om}\equiv
 \bigl\{\Om\times[0,1)\bigr\}\bigm\slash
 \bigl\{\Om\times\{0\}\bigr\}
\end{equation}
of a neighborhood $U\subset M$ of $\a$ onto a cone $K_{\Om}$ with smooth compact base $\Om$;
moreover, the mapping $\ph$ takes $\a$ to the vertex $\wt\a$ of the cone, and the restriction of
$\ph$ to the punctured neighborhood $\circU=U\setminus\a$ of the singular point is a
diffeomorphism.
\end{definition}

\medskip

\subsubsection*{Conical and cylindrical coordinates}
The coordinate on the interval $[0,1)$ in the representation~\eqref{zwezda} will be denoted
by~$r$. Points of $\Om$ (and sometimes local coordinates on $\Om$) will be denoted by $\om$. The
coordinates $(r,\om)$ are called the \textit{conical coordinates}. However, we mainly use the
\textit{cylindrical} coordinates related to the conical coordinates by the change of variables
$r=e^{-t}$. In the cylindrical representation, the smooth open manifold $\circM$ looks like a
manifold with a cylindrical end, and the punctured neighborhood $\circU$ is represented as the
direct product
\begin{equation}\label{pryapro}
 \circU=\Om\times\RR_+.
\end{equation}
The function $t:\circU\lra\RR_+$ well defined on $\circU$ by the direct product structure can be
continued to a smooth function $t:\circM\lra\RR$ everywhere nonpositive outside $\circU$. We
assume that the continuation is chosen and fixed.

\medskip

\subsubsection*{Weighted Sobolev spaces}
On $M$ we take some Riemannian metric $dx^2$ that has the form~\eqref{pryapro}
\begin{equation}\label{metric1}
dx^2=dt^2+d\om^2
\end{equation}
in the decomposition~\eqref{pryapro}, where $d\om^2$ is a Riemannian metric on $\Om$. The
\textit{weighted Sobolev space} $H^{s,\ga}(M)$, $s,\ga\in\RR$, is defined as the completion of
the space $C_0^\infty(\circM)$ of smooth compactly supported functions on $\circM$ with respect
to the norm
$$
\norm{u}_{s,\ga}=\left\{\int \left|(1-\triangle_M)^{s/2}e^{\ga t}u\right|^2\,d\mu\right\}^{1/2},
$$
where $\triangle_M$ and $d\mu$ are, respectively, the Beltrami--Laplace operator and the measure
associated with the metric~\eqref{metric1}. The numbers $s$ and $\ga$ are called the
\textit{smoothness exponent} and the \textit{weight exponent}, respectively. We deal only with
pseudodifferential operators in the Sobolev spaces $H^{s}(M)\equiv H^{s,0}(M)$ with zero weight
exponent; the general case can be reduced to this with the help of the isomorphism\footnote{Note
that the conjugation by this isomorphism provides an equivalence of pseudodifferential
quantizations corresponding to the weight exponents $\gamma$ and $0$.}
\begin{equation*}
 e^{-\ga t}:H^s(M)\lra H^{s,\ga}(M).
\end{equation*}
Moreover, we mainly deal with zero-order pseudodifferential operators, and so our main space is
$L^2(M)\equiv L^2(M,d\mu)=H^0(M)$. Occasionally, we shall use other spaces from the scale
$\{H^s(M)\}$.

\medskip

\subsubsection*{The cosphere bundle}
The cosphere bundle $S^*M$ of the manifold $M$ with a conical
singular point is of interest to us as the space on which the
interior principal symbols of pseudodifferential operators are
defined. This space is defined in the standard manner as the
quotient of the complement to the zero section in the
\textit{compressed cotangent bundle} $T^*M$ by the action of the
multiplicative group $\RR_+$ of positive numbers. An intrinsic
geometric definition of the compressed cotangent bundle can be
found, say, in Melrose~\cite{Mel4}. We do not reproduce it here
but just describe $S^*M$ in cylindrical coordinates, which is
sufficient for our aims. Namely, we consider the cosphere bundle
$S^*\circM$ and compactify it by attaching the manifold
$$
\bigl[(T^*\Om\times\RR)\setminus\{0\}\bigr]/\RR_+
$$
(where $\{0\}$ is the zero section of the vector bundle $T^*\Om\times\RR\lra\Om$) at $t=\infty$
as follows. The direct product structure on $\circU$ induces a canonical isomorphism
$$
S^*_{(\om,t)}M\simeq\bigl[(T^*_\om\Om\times\RR)\setminus\{0\}\bigr]/\RR_+
$$
for each point $(\om,t)\in\circU$. Now we say that a sequence $(x_k,\xi_k)\in S^*\circM$, where
$x_k=(\om_k,t_k)$ and $\xi_k\in\bigl[(T^*_{\om_k}\Om\times\RR)\setminus\{0\}\bigr]/\RR_+$, tends
to $(\om,\xi)\in\bigl[(T^*\Om\times\RR)\setminus\{0\}\bigr]/\RR_+$ if $t_k\lra\infty$,
$\om_k\lra\om$, and $\xi_k\lra\xi$. The compactified space $S^*M$ thus obtained is a smooth
manifold with boundary
\begin{equation}\label{passtar}
\pa S^*M= \bigl[(T^*\Om\times\RR)\setminus\{0\}\bigr]/\RR_+.
\end{equation}
(The smooth structure on $S^*M$ is defined via the radial coordinates.)

\medskip

\subsubsection*{The double}

In the subsequent argument, it is now and then convenient to
interpret operators on $M$ whose integral kernels are compactly
supported in $\circM\times\circM$ as operators on some closed
manifold. The simplest way to obtain this closed manifold is to
cut away the cylindrical end (say, at $t=10$) and attach the
second copy of the resulting manifold with boundary.  This
procedure gives a closed compact manifold, which is called the
\textit{double} of $M$ and denoted by $2M$. (Note that in contrast
with the papers~\cite{NaSt3,ScSS11}, dealing with index problems
for operators with symmetry conditions, where the passage to the
double has a topological meaning, in our case this passage is a
purely technical tool of analysis.) We can perform a similar
cut-and-paste procedure with $S^*M$, which results in $S^*2M$.
(Note, however, that the fibres of the cosphere bundle must be
glued together via the involution $p\longmapsto-p$ rather than
identically, where $p$ is the dual variable of $t$. This is
because the directions of the $t$-axis on the first and second
copies are opposite.)

\medskip

\subsubsection*{Pseudodifferential operators and symbols}
We deal with scalar classical pseudodifferential operators of order zero. The algebra of such
operators will be denoted by $\op{\Psi}(M)$. Every pseudodifferential operator
\begin{equation*}
 \wh{a}:L^2(M)\lra L^2(M)
\end{equation*}
is uniquely determined modulo compact operators by its \textit{principal symbol}, which is the
pair
\begin{equation*}
 \si(\wh{a})=\bigl(A(x,\xi),\bA(p)\bigr)
\end{equation*}
consisting of the \textit{interior symbol} $A(x,\xi)$, which is a $C^\infty$ function on $S^*M$,
and the \textit{conormal symbol} $\bA(p)$, which is a family of classical zero-order
pseudodifferential operators with parameter $p\in\RR$ on the manifold $\Om$ in the sense of
Agranovich--Vishik~\cite{AgVi1}. (In particular, this means that the operator
\begin{equation*}
 \wDo:L^2(C_\Om)\lra L^2(C_\Om)
\end{equation*}
on the cylinder $C_\Om=\Om\times\RR$ is a translation-invariant classical zero-order
pseudodifferential operator. This operator will sometimes be also referred to as the conormal
symbol.) The elements of this pair must satisfy the matching condition
\begin{equation*}
 \si(\bA(p))=A(x,\xi)\big|_{\pa S^*M}\,,
\end{equation*}
where $\si(\bA(p))$ is the restriction to the spheres of the principal symbol of $\bA(p)$ treated
as an operator with parameter.

\medskip

\subsubsection*{The Calkin algebra and ellipticity}
The principal symbols $\si(\wh{a})$ thus defined are just elements of the Calkin algebra
\begin{equation*}
 \cA_\infty=\op{\Psi}(M)\bigm\slash\{\cK\cap\op{\Psi}(M)\}
\end{equation*}
of the algebra of zero-order classical pseudodifferential operators in $L^2(M)$. (Here $\cK$ is
the ideal of compact operators in $L^2(M)$.) The algebra $\cA_\infty$ is a dense local subalgebra
of the symbol $C^*$-algebra $\cA$ obtained as the closure of $\cA_\infty$ in the direct sum
$C(S^*M)\oplus \cL(L^2(C_\Om))$, where $\cL(L^2(C_\Om))$ is the algebra of continuous linear
operators in $L^2(C_\Om)$. The multiplication law in $\cA$ is defined separately on each
component; it is given by pointwise multiplication of functions on the first component (interior
symbols) and by multiplication of operator families on the second component. The
\textit{pseudodifferential quantization}
\begin{equation*}
 \tau_{\op{\Psi}}:\cA_\infty\lra \op{\Psi}(M)
\end{equation*}
on $M$ is a continuous linear mapping such that
\begin{equation*}
\si\circ \tau_{\op{\Psi}}=\op{id}\,.
\end{equation*}
The construction of this mapping, which uses
partitions of unity, can be found in the cited
literature. This mapping is unique and
multiplicatively homomorphic modulo compact operators.

The invertibility of the principal symbol (that is, the nonvanishing of $A(x,\xi)$ and the
invertibility of $\bA(p)$ for each $p\in\RR$) is a sufficient and necessary condition for the
operator $\wh{a}$ to be Fredholm. Invertible elements in $\cA$ (and in $\op{M}(\cA)$ are also
said to be \textit{elliptic}.

\subsection{A Toeplitz quantization of $\cA$}

Let us now present our results.

To specify a Toeplitz quantization of the algebra $\cA$, one should describe a Hilbert
$\cA$-module $\gW$ and an orthogonal projection $P$ in $\gW$ almost commuting with the action of
the algebra. Let us construct both objects.

\medskip

\subsubsection*{The $\cA$-module $\wt{\cH}$}

To choose this module, we use the following motivation.
Irreducible representations of $\cA$ fall into two series -- see
Plamenevsky-Senichkin~\cite{PlSe2,PlSe1}. One series consists of
one-dimensional representations parametrized by points $(x,\xi)\in
S^*\circM$ and given by the formula
\begin{equation*}
 \mu_{(x,\xi)}(a)=A(x,\xi),\quad a=(A(x,\xi),\bA(p))\in\cA\,.
\end{equation*}
The other series consists of the infinite-dimensional representations
\begin{equation*}
 \mu_p(a)=\bA(p),\quad p\in\mathbb{R}
\end{equation*}
in the Hilbert space $L^2(\Om)$. This suggests that the Hilbert
$\cA$-module $\gW$ (see~\eqref{toep}) occurring in the definition
of the Toeplitz quantization can be obtained as a direct integral
of these representations. More precisely, we take
\begin{equation*}
 \wt{\cH}=L^2(S^*\circM)\oplus L^2(C_\Om)
\end{equation*}
(where the measure on $S^*\circM$ is translation invariant on the cylindrical end) and define the
action of $\cA$ on $\gW$ componentwise by the formula
\begin{equation*}
 a(u\oplus v)=A(x,\xi)u\oplus\wAo v,\quad a=(A(x,\xi),\bA(p))\in\cA\,.
\end{equation*}
The $\cA$-module thus defined is faithful.

\medskip

\subsubsection*{The subalgebra $\wt\cA$}
We construct a Toeplitz quantization (and prove its equivalence to the pseudodifferential
quantization) not for the entire algebra $\cA$ but for the subalgebra $\wt\cA\subset\cA$
consisting of symbols $a=(A(x,\xi),\bA(p))$ stabilizing at infinity, i.e., such that
$A(x,\xi)\equiv A(\om,t,\xi)$ is independent of $t$ for $t\ge0$. The passage to this subalgebra
is sufficient at least in index theory, since $\wt\cA$ has the following properties:
\begin{itemize}
  \item an element is invertible in $\op{M}(\wt\cA)$ if and only if
        it is invertible in $\op{M}(\cA)$;
  \item each elliptic element $a\in\op{M}(\cA)$ is homotopic
        in $\op{M}(\cA)$ via elliptic elements with the same
        conormal symbol to an element of $\op{M}(\wt\cA)$.
\end{itemize}
The homotopy mentioned in the latter item ``sweeps'' all the variation of $A(\om,t,\xi)$ in $t$
on the cylindrical end into a small half-neighborhood ($-\e<t<0$) of the section $\{t=0\}$ of
$M$.

\medskip

\subsubsection*{The Toeplitz projection $P$  and the pseudo-Guillemin transform}

In this subsection we construct the projection $P$ by describing its range ${\cH}\subset
\wt{\cH}$ as the range of a continuous mapping
\begin{equation*}
 \Ga:L^2(M)\lra \wt{\cH},
\end{equation*}
which proves to be an almost isomorphism of $L^2(M)$ onto ${\cH}$. This mapping will be called
the \textit{pseudo-Guillemin transform} (or the Guillemin transform for manifolds with
singularities), since, as we shall see shortly, it is a counterpart of the ordinary Guillemin
transform~\cite{Bout3,BoGu1,Gui1} for smooth compact manifolds and specifies an equivalence
between the Toeplitz quantization corresponding to $P$ and the pseudodifferential quantization.

Thus let us describe $\Ga$. On the real line $\RR$, we consider a partition of unity
\begin{equation*}
 1=\bigl(\chhi1(t)\bigr)^2+\bigl(\chhi2(t)\bigr)^2,
\end{equation*}
such that the $\chhi{j}(t)$ are smooth real-valued functions and
\begin{equation*}
 \chhi1(t)=
 \begin{cases}
 0,& t\le1,\\
 1,& t\ge3.
 \end{cases}
\end{equation*}
Next, let $\psi(t)$ be a real-valued function such that $\supp\psi\subset\{t<4\}$ and
$\psi(t)\chhi2(t)=\chhi2(t)$. All these functions can be viewed as functions on $M$ and $C_\Om$
as well as on the double $2M$ if we extend them as continuous functions beyond the cut $t=10$ by
constant values on the newly attached second copy of $M$. Using the natural projections, we also
lift these functions to $S^*M$, $S^*C_\Om$, and $S^*2M$.

Let
\begin{equation}\label{ordgui}
 \Tau :L^2(2M)\lra L^2(S^*2M)
\end{equation}
be the Guillemin transform~\cite{Gui1} for the compact closed $C^\infty$ manifold $2M$. We define
a mapping
\begin{equation}\label{pseudoG0}
 \Ga:L^2(M)\lra \wt{\cH}=L^2(S^*M)\oplus L^2(C_\Om)
\end{equation}
by the formula
\begin{equation}\label{pseudoG}
  \Ga\ph=\psi \Tau \chhi2\ph\oplus\chhi1\ph,\quad\ph\in L^2(M).
\end{equation}
Here the right-hand side is well defined. Indeed, the function $\chhi2\ph$ is supported in
$\{t<3\}$ and hence can be treated as a function on $2M$, whereby we can apply the Guillemin
transform $\Tau $. Next, the multiplication by $\psi$ permits us to understand $\psi \Tau
\chhi2\ph$ as a function on $S^*M$ (supported in $\{t<4\}$) rather than $S^*2M$. Likewise,
$\chhi1\ph$ can be viewed as a function on the cylinder $C_\Om$.

\begin{definition}
The continuous mapping~\eqref{pseudoG0}, \eqref{pseudoG} is called the \textit{pseudo-Guillemin
transform} on a manifold $M$ with singularities.
\end{definition}

Now we state our main theorem.

\begin{theorem}\label{main-theorem}
The following assertions hold.
\begin{itemize}
  \item[\rom1.] The range ${\cH}=\im(\Ga)$ of the
                mapping $\Ga$ is closed.
  \item[\rom2.] The mapping $\Ga:L^2(M)\lra {\cH}$
                is Fredholm and almost unitary.\footnote{The latter assertion means that
$$
\Ga^*\Ga=\id_{L^2(M)}+K_1,\quad \Ga\Ga^*=\id_{\cH}+K_2,
$$
             where $K_1$ and $K_2$ are compact operators.}
  \item[\rom3.] The orthogonal projection $P$ on ${\cH}$
  in $\wt{\cH}$ satisfies the
  condition
\begin{equation*}
   [P,a]\in\mathfrak{S}_k\quad\text{for $k>2n-1$ and
  for arbitrary $a\in\wt\cA$.}
\end{equation*}
  Thus,
  $P$ determines a $k$-summable Toeplitz quantization of $\wt\cA$.
  \item[\rom4.] One has
$$
\Ga^*a\Ga=\tau_{\op{\Psi}}(a)
$$
modulo compact operators of the von Neumann--Schatten class $\mathfrak{S}_k$, so that $\Ga$
determines an equivalence between the Toeplitz quantization associated with $P$ and the
pseudodifferential quantization $\tau_{\op{\Psi}}$.
\end{itemize}
\end{theorem}

Thus the theorem states that the pseudo-Guillemin transform $\Ga$ defines a representation of the
pseudodifferential quantization on a manifold with isolated singularities as a Toeplitz Fredholm
module of order $k>2n-1$.

\subsection{Proof of Theorem~\ref{main-theorem}}

First we recall the main properties of the usual Guillemin transform, which will be useful in our
proof.

\begin{proposition}[see~\cite{Gui1}]
The Guillemin transform~\eqref{ordgui} has the following properties\rom:
\begin{itemize}
  \item[\rom{(1)}] $\Tau $ is an operator of order zero in Sobolev scales,
  that is, $\Tau :H^s(2M)\lra H^s(S^*2M)$ is bounded for all $s$\rom;
  \item[\rom{(2)}] the range $\im(\Tau )$ of $\Tau $ is closed\rom;
  \item[\rom{(3)}] $\Tau ^*\Tau =\id$, and $\Tau \Tau ^*=\Pi$ is
  the projection onto $\im(\Tau )$\rom;
  \item[\rom{(4)}]\label{4} for every smooth function $b$ on $S^*2M$, the commutator\footnote{%
Here and in the following, $b$ is interpreted as
        the operator of point-wise multiplication by $b$.}
                   $[\Pi,b]$ is an  order $-1$ operator in the Sobolev scale
        on $S^*2M$ and belongs to the von Neumann--Schatten class
                   $\mathfrak{S}_k(L^2(S^*2M))$ for $k>2n-1$\rom;
  \item[\rom{(5)}]\label{5} for every smooth function $b$ on $S^*2M$ one has
\begin{equation*}
 \Tau ^*b\Tau =\wb+K,
\end{equation*}
where $\wb$ is an arbitrary pseudodifferential operator of order
zero on $2M$ with principal symbol $b$ and $K$ is a
pseudodifferential operator of order $-1$ in the Sobolev scale on
$2M$  thus belonging the von Neumann--Schatten class
$\mathfrak{S}_k(L^2(2M))$ for $k>n$.
\end{itemize}
\end{proposition}

Now we proceed directly to the proof of the main theorem. One has
$$
\Ga^*\Ga=\chhi2 \Tau ^*\psi^2\Tau \chhi2+\chhi1^2=1+K_1,
$$
where $K_1$ is a self-adjoint pseudodifferential operator of order
$-1$ on $M$ with integral kernel supported in the compact set
$\{t\le3\}\times\{t\le3\}$. It follows that $\Ga^*\Ga$ is Fredholm
and $\Ker \Ga^*\Ga=\Ker\Ga$ is finite-dimensional and consists of
functions supported in $\{t\le3\}$. By standard argument, we find
that $\im(\Ga)$ is closed.  Let $Q$ be the orthogonal projection
onto $\Ker\Ga$. Then the operator $B=\Ga^*\Ga+Q$ is an invertible
pseudodifferential operator of order $0$ and its inverse $B^{-1}$
has the form $B^{-1}=1+K_2$ with $K_2$ a self-adjoint
pseudodifferential operator of order $-1$ on $M$ whose integral
kernel is supported in the compact set $\{t\le3\}\times\{t\le3\}$.
(Hence both $K_1$ and $K_2$ belong to $\mathfrak{S}_k$.) Now the
projection on $\im(\Ga)$ has the form
$$
P=\Ga B^{-1}\Ga^*=\begin{pmatrix}
  \psi \Tau \chhi2B^{-1}\chhi2 \Tau ^*\psi
     &
   \psi \Tau \chhi2B^{-1}\chhi1 \\
  \chhi1B^{-1}\chhi2\Tau ^*\psi
     &
  \chhi1 B^{-1}\chhi1
                  \end{pmatrix}.
$$
Next,
\begin{align*}
\{\Ga\Ga^*\}|_{\cH}&=\{\Ga B^{-1}\Ga^*\Ga\Ga^*\}|_{\cH}\\
       &=\{\Ga B^{-1}(1+K_1)\Ga^*\}|_{\cH}\\
       &=(1+K)|_{\cH},
\end{align*}
where $K=\Ga B^{-1}K_1\Ga^*$ is compact, and we have proved
assertions~1 and~2. Since $B^{-1}$ differs from the identity
operator by an operator in  $\mathfrak{S}_k$, we can safely
replace the operator $P$ by
$$
\wP=\Ga \Ga^*=\begin{pmatrix}
  \psi \Tau \chhi2^2 \Tau ^*\psi
     &
   \psi \Tau \chhi2\chhi1 \\
  \chhi1\chhi2\Tau ^*\psi
     &
  \chhi1^2
                  \end{pmatrix}
$$
in the proof of assertion~3 and estimate the commutator $[\wP,a]$.
Moreover, owing to the presence of the cutoff function $\psi$ we
can assume that the interior symbol $A$ vanishes for $t>9$ and
hence $\wP$ can be treated as an operator acting in
$L^2(S^*2M)\oplus L^2(C_\Om)$ when we proceed with the estimates.
We have
$$
[\wP,a]=\begin{pmatrix}
  [\psi \Tau \chhi2^2 \Tau ^*\psi,A]
     &
   \psi \Tau \chhi2\chhi1\bhA -A\psi \Tau \chhi2\chhi1 \\
  \chhi1\chhi2\Tau ^*A-\bhA \chhi1\chhi2\Tau ^*\psi
     &
  [\chhi1^2,\bhA ]
                  \end{pmatrix},
$$
where we for brevity write
$$
\bhA =\wAo.
$$
Now, by the properties of $\Tau $,
$$
[\psi \Tau \chhi2^2 \Tau ^*\psi,A]\equiv[\Pi,A]\equiv0 \,.
$$
(Here and in the following $\equiv$ stands for equality modulo
operators in  $\mathfrak{S}_k$.) Next, we shall estimate the
commutator\footnote{Similar estimates of commutators of
pseudodifferential operators and cutoff functions on the infinite
cylinder can be found in~\cite{PlRo4}.} $[\chhi1^2,\bhA ]$. To
this end, we introduce a smooth partition of unity
$$
\ph_0(t)+\ph_1(t)+\ph_2(t)=1
$$
on the real line such that the following properties hold: (a) $\chhi1^2\ph_1=\ph_1$; (b)
$\chhi1^2\ph_2=0$; (c) $\supp \ph_0$ is contained in the interval $[0,10]$ and $\ph_0=1$ on
$\supp \chhi1'$).

Now we represent $\bhA $ in the form
$$
\bhA =\sum\limits_{j,k=0}^2\ph_j\bhA \ph_k
$$
and compute the commutator as follows:
\begin{equation*}
\begin{split}
[\chhi1^2,\bhA ]&=\sum\limits_{j,k=0,1,2}[\chhi1^2,\ph_j\bhA \ph_k]\\
&=[\chhi1^2,\ph_0\bhA \ph_0]-[\chhi2^2,\ph_0\bhA \ph_1]
 +[\chhi1^2,\ph_0\bhA \ph_2]\\
 &\quad{}-[\chhi2^2,\ph_1\bhA \ph_0]+[\chhi1^2,\ph_2\bhA \ph_0]
 +\sum\limits_{j,k=1,2}(-1)^{j}[\chhi{3-j}^2,\ph_j\bhA \ph_k].
\end{split}
\end{equation*}
Here we have used the fact that $\chhi1^2+\chhi2^2=1$.

Now the first term can be interpreted as an order $-1$
pseudodifferential operator  on $2M$ (the support of its Schwartz
kernel with respect to the variable $t$ is contained in the square
$[0,10]\times[0,10]$) and hence belongs to the desired von
Neumann--Schatten class $\mathfrak{S}_k$. All other terms have the
form $a\bhA b$, where $a=a(t)$ and $b=b(t)$ are smooth functions
with disjoint supports,  constant at infinity. For example,
$$
[\chhi2^2,\ph_0\bhA \ph_1]=\chhi2^2\ph_0\bhA \ph_1
 -\ph_0\bhA \ph_1\chhi2^2=\chhi2^2\ph_0\bhA \ph_1;
$$
here we have $a=\chhi2^2\ph_0$ and $b=\ph_1$.

The integral kernel of such an operator has the form
$$
K(x,x')=a(t)b(t')K_0(\om,\om',t,t')\,,
$$
where $K_0(\om,\om',t,t')$ is the integral kernel of $\bhA $. The kernel $K_0(\om,\om',t,t')$ is
smooth outside the diagonal $t=t'$ and decays more rapidly than an arbitrary power of
$|t-t'|^{-1}$. With regard for the arrangement of supports of $a(t)$ and $b(t)$, we have
$$
1+|t|+|t'|\le C|t-t'|
$$
with some constant $C$ on
$$
\supp K(x,x')\subset \supp a(t)\times \supp
b(t')\times\Om\times\Om\,.
$$
It follows that $K(x,x')$ is everywhere smooth and satisfies the estimates
$$
|K(x,x')|\le C_N(1+|t|+|t'|)^{-N}
$$
for all $N$; similar estimates are valid for the derivatives of $K(x,x')$. We conclude that the
operator $a\bhA b$ belongs to all von Neumann--Schatten classes. Thus, we have obtained the
desired estimate for the commutator $[\chhi1^2,\bhA ]$.

Consider the operator $\psi \Tau \chhi2\chhi1\bhA -A\psi \Tau \chhi2\chhi1$. We have
$$
 \psi \Tau \chhi2\chhi1\bhA -A\psi \Tau \chhi2\chhi1\equiv
 \psi \Tau \chhi2\chhi1\bhA -\psi \Tau \wA\chhi2\chhi1\,,
$$
where $\wA$ is some pseudodifferential operator with symbol $A$,
and the desired assertion follows (since $\si(\bhA)=A$ on
$\supp\chhi2\chhi1$) by an argument similar to the preceding one.
The estimate for the lower left entry of the commutator is
similar.

It remains to prove assertion~4. We have
$$
\Ga^*a\Ga=\chhi2\Tau ^*\psi A\psi \Tau \chhi2+\chhi1\bhA \chhi1 \equiv
 \chhi2^2\wA+\chhi1^2\bhA \equiv\tau_\Psi(a)\,,
$$
as follows from our conditions on the symbols in $\wt\cA$. The
proof is complete.\qed

\section{Equivalent Toeplitz quantizations}

\subsection{The construction via a bounded resolution}

We start with an important remark. In what follows, we sometimes
use a slight generalization of the construction of Toeplitz
quantization~\eqref{toep}. Suppose that the space $\wt\cH$ where
the projection $P$ defining the quantization acts is itself a
quantization space rather than a Hilbert $\cA$-module, i.e., is
equipped with an \textit{almost} representation\footnote{A linear
mapping $\tau':\cA\lra B(\wt\cH)$ is called an almost
representation if it is an algebra homomorphism modulo compact
operators: $\tau'(ab)=\tau'(a)\tau'(b)+K$, where $K$ is a compact
operator.} (quantization) $\tau'$ of the algebra $\cA$. (This
quantization may well be Toeplitz itself.) Let $P$ almost commute
with the quantization, i.e., satisfy
\begin{equation*}
    P\wt\tau'(a)\equiv\wt\tau'(a)P,\quad a\in\cA\,,
\end{equation*}
where $\equiv$ stands for equality modulo compact operators. (This condition is independent of
the specific choice of the representative $\wt\tau'$.) Then a slight modification of
formula~\eqref{toep} permits one to define a Toeplitz quantization $\tau$ in $\cH$ by setting
\begin{equation}\label{toep-1}
\begin{CD}
    \wt\tau(a)=P\wt\tau'(a)\,.
\end{CD}
\end{equation}
The quantization $\tau$ is well defined (since it is independent of the choice of a
representative $\wt\tau'$ of $\wt\tau$) and will be denoted by $P\tau'$.
If $\tau'$ is a Toeplitz
quantization associated with a projection $P'$, then $P\tau$ is associated with $PP'$.

\medskip

\subsubsection{Admissible operators}

Now let us describe the class of operators that will be used in our resolutions.

\begin{definition}\label{admi}
Let $\cH_1$ and $\cH_2$ be Hilbert $\cA$-modules. A bounded linear
operator
\begin{equation*}
    D:\cH_1\lra \cH_2
\end{equation*}
is said to be \textit{admissible} if the following conditions hold.
\begin{enumerate}
    \item The range of $D$ is closed.
    \item The operator $D$ almost commutes with the action of the algebra,
    i.e.,
\begin{equation*}
    Da-aD\in \cK(\cH_1,\cH_2),\quad a\in\cA\,.
\end{equation*}
\end{enumerate}
\end{definition}

In a similar way, we define admissible operators for
the case in which $\cH_1$ and $\cH_2$ are not
$\cA$-modules but quantizations of the algebra $\cA$
acting in these spaces.

Note that a projection is admissible if and only if it
almost commutes with the action (or quantization) of
the algebra $\cA$, i.e., specifies a Toeplitz
quantization.

\begin{proposition}\label{odin}
Let
\begin{equation*}
    D:\cH_1\lra \cH_2
\end{equation*}
be an admissible operator. Then the following assertions hold.
\begin{enumerate}
    \item[\rom{1)}]
    The operators $D^*$, $D^*D$, and $DD^*$, as well as the projections
    on the kernel, cokernel \rom(the kernel of the adjoint operator\rom),
    range, and corange of $D$ are admissible operators.
    \item[\rom{2)}] If
\begin{equation*}
 P_1:\cH_1\lra\cH_1,\quad P_2:\cH_2\lra\cH_2
\end{equation*}
    is a pair of admissible projections commuting with $D$,
    i.e.,
\begin{equation*}
    DP_1-P_2D=0,
\end{equation*}
    then the operator
\begin{equation*}
    P_2DP_1:\im(P_1)\lra \im(P_2)
\end{equation*}
    \rom(which is the restriction of
    $D$ to the range of $P_1$\rom)
    is also an admissible operator.
    \item[\rom{3)}] If, moreover, $D$ is self-adjoint, then its positive
    spectral projection $P_+(D)$ is an admissible operator.
\end{enumerate}
\end{proposition}

\begin{proof}
1) The fact that the ranges of $D$ and $D^*$ are or are not closed simultaneously is well known
in functional analysis. Under our assumptions, the ranges of $D^*D$ and $DD^*$ coincide with the
ranges of $D^*$ and $D$, respectively,  and hence are also closed. The ranges of bounded
projections are always closed. The assertion that these operators commute with the action of
$\cA$ modulo compact operators is nontrivial only for the projections. It suffices to consider
the projections on $\Ker D$ and $\Coker D$. Consider the first projection. (The second projection
can be treated in a similar way.) The kernel of $D$ coincides with the kernel of the self-adjoint
positive semidefinite operator $D^*D$. Since the range of $D^*D$ is closed, it follows that the
restriction of $D^*D$ to the orthogonal complement of $\Ker D$ is boundedly invertible by
Banach's theorem. Hence the operator $z-D^*D$ is invertible for nonzero $z\in\CC$ sufficiently
small in absolute value. Thus, the projection $P$ on $\Ker D$ can be expressed by the Cauchy
integral formula
\begin{equation*}
    P=\frac1{2\pi i} \oint\frac{dz}{z-D^*D},
\end{equation*}
where the integration is carried out over a circle of small radius centered at zero. Accordingly,
the commutator of $P$ with elements of $\cA$ is expressed by the formula
\begin{align*}
    [a,P]&=\frac1{2\pi i} \oint [a,(z-D^*D)^{-1}]\,dz \\
         &=\frac1{2\pi i} \oint (z-D^*D)^{-1}[a,D^*D](z-D^*D)^{-1}\,dz ,
\end{align*}
which readily implies the desired compactness.

2) Under our assumptions, the operator $P_2DP_1$ almost commutes with the action of $\cA$, and
the fact that its range is closed is obvious.

3) Since the range of $D$ is closed, it follows that zero is either in the resolvent set of $D$
or an isolated point of spectrum of $D$. Consequently, the positive spectral projection can be
expressed by a Cauchy type integral and we can argue as in the proof of~1).
\end{proof}

\medskip

\subsubsection*{Resolution and the associated self-adjoint operator}

Let $$P:\cH_0\lra\cH_0$$ be an admissible projection in a Hilbert $\cA$-module $\cH_0$. The
corresponding Toeplitz quantization will be denoted by $\tau$.
\begin{definition}
An \textit{admissible resolution} of length $n$ of the projection $P$ is an exact sequence of
Hilbert spaces
\begin{equation}\label{reso-1}
    \begin{CD}
    0@>>>\im(P)@>>>\cH_0@>A_0>>\cH_1@>A_1>>\dotsc@>A_{n-1}>>\cH_n\,,
    \end{CD}
\end{equation}
where all spaces are Hilbert $\cA$-modules and all mappings are admissible operators.
\end{definition}

A trivial example is the resolution of length $0$
\begin{equation}\label{reso-2}
    \begin{CD}
    0@>>>\im(P)@>>>\cH_0\,.
    \end{CD}
\end{equation}

For a given resolution~\eqref{reso-1}, we construct a self-adjoint operator $D$ in the direct sum
\begin{equation}\label{summan}
    \cH=\bigoplus_{k=0}^n\cH_n
\end{equation}
of Hilbert spaces as follows. First, we augment the exact sequence~\eqref{reso-1} with the
cokernel of $A_{n-1}$:
\begin{equation}\label{reso-3}
    \begin{CD}
    0@>>>\im(P)@>>>\cH_0@>A_0>>\cH_1@>A_1>>\dotsc\\
    &&&&@>A_{n-1}>>\cH_n
    @>>>\Coker A_{n-1}@>>>0\,.
    \end{CD}
\end{equation}
The projection in $\cH_n$ on the cokernel $\Coker A_{n-1}$ will be denoted by $\wt P$ and called
the projection \textit{dual} to $P$. By Proposition~\ref{odin}, this is an admissible projection.
Next, let $A$ be the differential of degree $+1$ defined in the graded sum~\eqref{summan} as the
direct sum of all operators $A_j$, and let $A^*$ be the adjoint operator (which is automatically
of degree $-1$). We set
\begin{equation}\label{dirac-abs}
    D=A+A^*+P+(-1)^{n+1}\wt P\,,
\end{equation}
where the operators $P$ and $\wt P$ act in the respective components $\cH_0$ and $\cH_n$ of the
direct sum. By the preceding reasoning, it is obvious that the operator $D$ is self-adjoint,
invertible, and admissible.

\medskip

\subsubsection*{Equivalence of quantizations}

If $\tau$ is a quantization of the algebra $\cA$, then by $[\tau]\in K^1(\cA)$ we denote its
equivalence class in the $K$-group of $\cA$ (see~\cite{Bla1}). Moreover, if $\tau$ is a Toeplitz
quantization determined by a projection $P$, then we freely write $[P]$ instead of $[\tau]$.

Let us now establish the main theorem of this section.

\begin{theorem}\label{34}
The positive spectral projection $P_+(D)$ of the operator $D$ is admissible. The Toeplitz
quantization corresponding to the projection $P_+(D)$ is equivalent to the Toeplitz quantization
defined by the projection $P$\rom:
\begin{equation*}
    [P]=[P_+(D)]\in K^1(\cA)\,.
\end{equation*}
\end{theorem}

\begin{proof}
The first assertion of the theorem follows from Proposition~\ref{odin},\,3). Let us prove the
equivalence of quantizations. We can break the resolution~\eqref{reso-1} at the $j$th space
$\cH_j$, thus obtaining an admissible resolution of length $j$. We denote the projection $\wt P$
obtained for the latter resolution by the same construction by $\wt P_j$, the direct
sum~\eqref{summan} with $n$ replaced by $j$ by $H_j$, and the corresponding operator $D$ by
$D_j$. We prove the theorem by induction on $j$.

(a) The inductive base. Let $j=0$. Then $H_0=\cH_0$, the augmented resolution~\eqref{reso-3} has
the form
\begin{equation*}
    \begin{CD}
    0@>>>\im(P)@>>>\cH_0@>>>\im(P)^{\perp}@>>>0\,,
    \end{CD}
\end{equation*}
the dual projection is $\wt P_0=1-P$, and
\begin{equation*}
    D_0=P-(1-P)=2P-1\,.
\end{equation*}
The positive spectral projection of the operator $D_0$ is just $P$, and so the assertion is
trivially true.

(b) The inductive step. Suppose that the assertion has already been proved for $j=k$, and let us
prove it for $j=k+1$. To this end, we write out the corresponding resolutions of lengths $k$ and
$k+1$ next to each other:
\begin{equation*}
    \begin{CD}
    0@>>>\im(P)@>>>\cH_0@>A_0>>\cH_1@>A_1>>\dotsc@>A_{k-1}>>\cH_k\,,\\
    0@>>>\im(P)@>>>\cH_0@>A_0>>\cH_1@>A_1>>\dotsc@>A_{k-1}>>\cH_k
    @>A_{k}>>\cH_{k+1}\,.
    \end{CD}
\end{equation*}
We expand the spaces $\cH_k$ and $\cH_{k+1}$ into the orthogonal sums
\begin{equation}\label{summy}
    \cH_k=R_k\oplus E_k\,,\quad
    \cH_{k+1}=R_{k+1}\oplus E_{k+1}\,,
\end{equation}
where $R_j=\im(A_{j-1})$, $j=k,k+1$. The terms in these expansions are not $\cA$-modules, but
they are equipped with Toeplitz quantizations determined in the ambient spaces $\cH_j$ by the
projections $\wt P_j$ (in the case of $E_j$) and $1-\wt P_j$ (in the case of $R_j$). One can
readily see that the operators $D_k$ and $D_{k+1}$ can be represented as the direct sums
\begin{equation}\label{pryam}
    D_k=B\oplus(\pm \wt P_k),\quad
    D_{k+1}=B\oplus
    \begin{pmatrix}
      0 & A^* & 0 \\
      A & 0 & 0 \\
      0 & 0 & \mp\wt P_{k+1}
    \end{pmatrix},
\end{equation}
where $B$ is an admissible self-adjoint operator acting in the direct sum $H_{k-1}\oplus R_k$,
$A$ and $A^*$ are the restrictions of the operators $A_k$ and $A_k^*$ to the subspaces $E_k\lra
R_{k+1}$ and $R_{k+1}\lra E_k$, respectively (these restrictions are still the adjoints of each
other), the signs $\pm$ and $\mp$ depend on the parity of $k$, and the projection $\wt P_j$ are
assumed to be restricted\footnote{Thus, we might well write the identity operators instead of
these projections in the formulas; we still write the projections so as to emphasize that the
$E_j$ are equipped with the corresponding Toeplitz quantizations.} to the subspaces $E_j$. Thus,
the positive spectral projections of these operators have the form
\begin{equation}\label{polozhit}
    P_+(D_k)=P_+(B)\oplus\wt P_k,\quad
    P_+(D_{k+1})=P_+(B)\oplus
    P_+\left[\begin{pmatrix}
      0 & A^* \\
      A & 0 \\
    \end{pmatrix}\right]
\end{equation}
for the upper signs and
\begin{equation}\label{otrizat}
    P_+(D_k)=P_+(B),\quad
    P_+(D_{k+1})=P_+(B)\oplus
    P_+\left[\begin{pmatrix}
      0 & A^* \\
      A & 0 \\
    \end{pmatrix}\right]\oplus\wt P_{k+1}
\end{equation}
for the lower signs.

We claim that the positive spectral projection of the self-adjoint operator
\begin{equation}\label{antidiag}
V=  \begin{pmatrix}
      0 & A^* \\
      A & 0 \\
    \end{pmatrix}:
    E_k\oplus R_{k+1}\lra E_k\oplus R_{k+1}
\end{equation}
defines a Toeplitz quantization equivalent to the quantization associated with the projection
$1-\wt P_{k+1}$. Indeed, first we note that the operator~\eqref{antidiag} is self-adjoint,
admissible and invertible. Next, consider the unitary transformation
\begin{equation}\label{preobr}
    U=\frac1{\sqrt2}\begin{pmatrix}
      1 & -(A^*A)^{-1/2}A^* \\
      (AA^*)^{-1/2}A & 1 \\
    \end{pmatrix}:
    E_k\oplus R_{k+1}\lra E_k\oplus R_{k+1}\,.
\end{equation}
It is admissible (by holomorphic functional calculus) and reduces $V$ to the block diagonal form
\begin{equation}\label{diag}
    UVU^{-1}=\begin{pmatrix}
      -(A^*A)^{1/2} & 0 \\
      0 & (AA^*)^{1/2} \\
    \end{pmatrix}.
\end{equation}
The positive spectral subspace of the operator~\eqref{diag} coincides with
\begin{equation*}
 R_{k+1}=\im(1-\wt P_{k+1})\,,
\end{equation*}
which gives the desired equivalence. By $\tau_k$ we denote the quantization
corresponding to the projection $P_+(D_k)$. In this notation, for the case of lower signs
in~\eqref{pryam}, we have
\begin{equation}\label{niz}
    [\tau_{k+1}]=[\tau_k]+[1-\wt P_{k+1}]+[\wt P_{k+1}]=[\tau_k]\,,
\end{equation}
since the elements $[1-\wt P_{k+1}]$ and $[\wt P_{k+1}]$ are the inverses of each other in the
$K$-group. For the case of upper signs, we have
\begin{equation}\label{verh}
    [\tau_{k+1}]=[\tau_k]+[1-\wt P_{k+1}]-[\wt P_{k}]\,.
\end{equation}
To complete the proof, let us show that
\begin{equation}\label{chu}
    [1-\wt P_{k+1}]=[\wt P_{k}]\,.
\end{equation}
Indeed, the operator $A$ is an admissible isomorphism
\begin{equation*}
    A:E_k\equiv \im(\wt P_{k})\lra R_{k+1}\equiv (1-\wt P_{k+1})\,,
\end{equation*}
whence~\eqref{chu} follows. The proof is complete.
\end{proof}

\subsection{The construction via an unbounded resolution}

Now let us consider the case in which one deals with an unbounded resolution of the projection
operator. Resolutions of such form naturally arise in applications.

\medskip

\subsubsection*{Main assumptions}

Consider the sequence
\begin{equation}\label{unbd-1}
    \begin{CD}
    0@>>>\cH_0@>A_0>>\cH_1@>A_1>>\dotsc@>A_{n-2}>>\cH_{n-1}@>A_{n-1}>>\cH_n@>>>0\,,
    \end{CD}
\end{equation}
where
\begin{enumerate}
    \item[1)] all $\cH_j$, $j=0,\dotsc,n$, are Hilbert $\cA$-modules;
    \item[2)] all $A_j$, $j=0,\dotsc,n-1$, are closed densely defined
    linear operators with closed range;
    \item[3)] the operators $A_j$ boundedly commute with the action of
    the algebra $\cA$, that is, for each $a\in\cA$ and each $j=0,\dotsc,n-1$
    one has $a \cD(A_j)\subseteq\cD(A_j)$ and the commutator $[a,A_j]$,
    well defined on $\cD(A_j)$ by virtue of this inclusion, is bounded
    and hence can be extended by closure to an everywhere defined bounded operator
    from $\cH_j$ to $\cH_{j+1}$, which will be denoted by the same symbol
    $[a,A_j]$.
\end{enumerate}
Further, we suppose that
\begin{enumerate}
    \item [4)] the sequence~\eqref{unbd-1} is a complex and, moreover, is
    exact in all but the extreme terms, that is,
\begin{equation}\label{exact}
    \Ker A_j=\im A_{j-1},\quad j=1,\dotsc,n-1.
\end{equation}
\end{enumerate}

By $P_0$ and $P_n$ we denote the orthogonal projections on $\Ker A_0$ in $\cH_0$ and $(\im
A_{n-1})^\perp$ in $\cH_{n}$, respectively; these projections correspond to the only nonzero
cohomology spaces of the complex~\eqref{unbd-1}. Just as in the bounded case, $P_n$ will be
referred to as the dual projection for $P_0$.

\medskip

Consider the Laplacians
\begin{equation}\label{lapl}
    \Delta_j=A_j^*A_j+A_{j-1}A_j^*:\cH_j\lra\cH_j,\quad
    \quad j=0,\dotsc,n\,,
\end{equation}
of the complex~\eqref{unbd-1}. Here, to unify the notation, we have set $A_{-1}=A_n=0$, so that
actually
\begin{equation*}
    \Delta_0=A_0^*A_0,\quad\Delta_n=A_{n-1}A_{n-1}^*.
\end{equation*}
The operators~\eqref{lapl} are densely defined and self-adjoint. Indeed, it follows from
condition~2) that $\im A_j^*$ is closed for all $j=0,\dotsc,n-1$, and with regard
to~\eqref{exact} one has the orthogonal expansions
\begin{equation}\label{expa}
    \cH_j=\im A_{j-1}\oplus\im A_j^*,\quad j=1,\dotsc,n-1
\end{equation}
(see~\cite[Theorem~IV.5.13]{Kat1}). One can readily see that for $j=1,\dotsc,n-1$ the Laplacian
$\Delta_j$ has the direct sum decomposition
\begin{equation}\label{dirsum}
    \Delta_j=A_{j-1}A_{j-1}^*\oplus A_j^*A_j
\end{equation}
corresponding to~\eqref{expa}. It remains to note that for a
closed densely defined operator $T$, the operator $T^*T$ is
densely defined and self-adjoint (see Kato,
\cite[Theorem~V.3.24]{Kat1}).

We impose the following condition on the Laplacians of the complex~\eqref{unbd-1}:
\begin{enumerate}
    \item[5)] for $j=1,\dotsc,n-1$, the Laplacian $\Delta_j$ is an
    operator with compact resolvent.
\end{enumerate}

\begin{remark}
The Laplacians $\Delta_0$ and $\Delta_n$ cannot have a compact resolvent except for the trivial
case in which the cohomology groups of the complex~\eqref{unbd-1} in dimensions $0$ and $n$ are
finite-dimensional and hence the corresponding projections $P_0$ and $P_n$ define trivial
quantizations. That is why we impose this condition only on the Laplacians in intermediate
dimensions. Condition~5) is empty for $n=1$. From now on, we assume that $n>1$; the subsequent
results remain valid for $n=1$ provided that condition~5) is replaced by the following condition:
\begin{enumerate}
    \item[5$'$)] the restrictions $\Delta_j|_{\im P_j}$, $j=0,1$, are
    operators with compact resolvent.
\end{enumerate}
\end{remark}

\medskip

\subsubsection*{Admissible resolutions}

Now we can introduce the following definition.

\begin{definition}
Suppose that conditions~1)--5) are satisfied. Then one says that the complex~\eqref{unbd-1} is an
\textit{admissible unbounded resolution} of length $n$ of the projection $P_0$.
\end{definition}

Just as in the bounded case, let $A$ be the differential of degree $+1$ defined in the graded sum
$\cH=\bigoplus_{k=0}^n\cH_n$ as the direct sum of all operators $A_j$, and let $A^*$ be the
adjoint operator. We define the rolled-up operator $D$ by the same formula
\begin{equation}\label{unb-dirac-abs}
    D=A+A^*+P_0+(-1)^{n+1} P_n
\end{equation}
as in the bounded case, where it is assumed that the operators $P$ and $\wt P$ act in the
respective components $\cH_0$ and $\cH_n$ of the direct sum and
\begin{equation*}
    \cD(D)=\cD(A)\cap\cD(A^*)\,.
\end{equation*}

Since the sequence~\eqref{unbd-1} is a complex, one can readily establish the following
assertion.

\begin{lemma}
The operator $D$ is self-adjoint.
\end{lemma}

\medskip
\begin{small}
\begin{proof}
Indeed, since $P_0$ and $P_n$ are self-adjoint and bounded, it suffices to prove the
self-adjointness of $B=A+A^*$ on the same domain $\cD(B)=\cD(A)\cap\cD(A^*)$. The range
\begin{equation*}
 \im A^*=\bigoplus_{j=0}^{n-1}\im A_j^*
\end{equation*} is closed, and so we have the orthogonal decomposition
\begin{equation*}
    \cH=\Ker A\oplus\im A^*.
\end{equation*}
Now $\Ker A\subset\cD(A)$ (tautologically) and $\im A^*\subset\cD(A^*)$ (since $(A^*)^2=A^2=0$),
and so we have the orthogonal decompositions
\begin{equation*}
    \cD(A)=\Ker A\oplus\cD,\quad \cD\subset\im A^*,\qquad
    \cD(A^*)=\wt\cD\oplus\im A^*,\quad\wt\cD\subset\Ker A\,,
\end{equation*}
where both embeddings are dense. It follows that
\begin{equation*}
    \cD(B)\equiv\cD(A+A^*)=\wt\cD\oplus\cD.
\end{equation*}

Now let us prove that $\cD(B^*)=\cD(B)$ (and then automatically $B=B^*$). The inclusion
$\cD(B)\subset\cD(B^*)$ is obvious, and we have only to prove the opposite inclusion. Let $u\in
D(B^*)$. This means by definition that there exists a constant $C$ such that
\begin{equation}\label{spitze}
    \lvert(u,Bv)\rvert\le C\norm{v}
\end{equation}
for each $v\in\cD(B)$. In particular, for arbitrary $v\in\wt\cD$ inequality~\eqref{spitze}
becomes
\begin{equation*}
    \lvert(u,A^*v)\rvert\le C\norm{v},
\end{equation*}
and now for arbitrary $x=v+w\in\cD(A^*)$, $w\in\im(A^*)$, we have $A^*x=A^*v$ and
\begin{equation*}
    \lvert(u,A^*x)=\lvert(u,A^*v)\rvert\le C\norm{v}\le C\norm{x}.
\end{equation*}
It follows that $u\in\cD(A^{**})=\cD(A)$. (Recall that $A$ is closed.) Next, for arbitrary
$v\in\cD$ inequality~\eqref{spitze} becomes
\begin{equation*}
    \lvert(u,Av)\rvert\le C\norm{v},
\end{equation*}
and now for arbitrary $x=w+v\in\cD(A)$, $w\in\Ker(A)$, we have
\begin{equation*}
    \lvert(u,Ax)=\lvert(u,Av)\rvert\le C\norm{v}\le C\norm{x}.
\end{equation*}
It follows that $u\in\cD(A^{*})$. We conclude that $u\in\cD(B)$, which completes the proof.
\end{proof}
\end{small}

\medskip

\subsubsection*{Main theorem}

Now we can state the main theorem.

\begin{theorem}\label{main}
Suppose that the sequence~\eqref{unbd-1} is an admissible unbounded resolution of length $n>1$ of
the projection $P_0$. Then the following assertions hold.
\begin{enumerate}
    \item[\rom{i)}] The projections $P_0$ and $P_n$, as well as the
    positive spectral projection $P_+(D)$ of the rolled-up operator $D$
    given by~\eqref{unb-dirac-abs}, are admissible.
    \item[\rom{ii)}] The Toeplitz quantizations defined by
    the projections $P_+(D)$ and $P$ are equivalent.
\end{enumerate}
\end{theorem}

\begin{proof}
We split the proof into two lemmas. For $j=0,\dotsc,n-1$, we define operators
\begin{equation}\label{bj}
    B_j:\cH_j\lra\cH_{j+1}
\end{equation}
by the formula
\begin{equation}\label{bj1}
    B_j=A_j(1+\Delta_j)^{-1/2}.
\end{equation}
\begin{lemma}\label{lemmaa}
The operators~\eqref{bj}, \eqref{bj1} are everywhere defined and bounded. Moreover, they are
admissible in the sense of Definition~\rom{\ref{admi}}, and
\begin{equation}\label{svoistva}
    \Ker B_j=\Ker A_j,\quad \im B_j=\im A_j.
\end{equation}
\end{lemma}
\begin{corollary}
The sequence
\begin{equation}\label{unbd-2}
    \begin{CD}
    0\lra\im P_0\lra\cH_0@>B_0>>\cH_1@>B_1>>\dotsc@>B_{n-2}>>\cH_{n-1}@>B_{n-1}>>\cH_n
    \end{CD}
\end{equation}
is a \rom(bounded\rom) admissible resolution of length $n$ of the projection $P_0$ with the same
dual projection $P_n$.
\end{corollary}
Indeed, this readily follows from~\eqref{svoistva}.
\begin{proof}[Proof of Lemma~\rom{\ref{lemmaa}}]
In intermediate dimensions, we have the decomposition \eqref{dirsum} corresponding to the direct
sum~\eqref{expa}, and accordingly,
\begin{equation}\label{dirsum1}
  (1+\Delta_j)^{-1/2}=(1+A_{j-1}A_{j-1}^*)^{-1/2}\oplus(1+ A_j^*A_j)^{-1/2}.
\end{equation}
(This readily follows from the properties of the functional calculus of self-adjoint operators;
e.g., see~\cite[\S\,VI.5.2]{Kat1} or \cite[Chap.~9]{RiNa1}.) Since $A_j$ is zero on the first
component of the expansion~\eqref{expa}, it follows that we can safely replace
$(1+A_{j-1}A_{j-1}^*)^{-1/2}$ by $1$ in the formula for $B_j$, so that we have
\begin{equation*}
    B_j=A_j(1+A_j^*A_j)^{-1/2}, \quad j=0,\dotsc,n-1.
\end{equation*}
(For $j=0$ this is valid automatically, since $\Delta_0=A_0^*A_0$.) Now we note that the
following assertion is true.

\textit{For a closed densely defined operator $T$, the range of the operator $(1+T^*T)^{-1/2}$
coincides with the domain of $T$.}

Indeed,
\begin{equation*}
 \im (1+T^*T)^{-1/2}=\cD\bl((1+T^*T)^{1/2}\br).
\end{equation*}
Next, the domain $\cD(1+T^*T)=\cD(T^*T)$ is a core of $(1+T^*T)^{1/2}$ as well as of
$T$~\cite[Chap.~V, Theorem~3.24 and Lemma~3.38]{Kat1}. Let us prove that
\begin{equation*}
 \cD\bl((1+T^*T)^{1/2}\br)=\cD(T).
\end{equation*}
If $u_k\in\cD(1+T^*T)$ is a sequence such that
\begin{equation*}
 u_k\to u\quad\text{and}\quad (1+T^*T)^{1/2}u_k\to v,
\end{equation*} then
\begin{equation*}
    \norm{(1+T^*T)^{1/2}(u_k-u_s)}^2=
    \norm{u_k-u_s}^2+\norm{Tu_k-Tu_s}^2\to0
\end{equation*}
as $n,m\to\infty$, and hence the sequence $Tu_k$ is convergent. Thus,
\begin{equation*}
 \cD\bl((1+T^*T)^{1/2}\br)\subset\cD(T).
\end{equation*}
Reversing the argument, we see that
\begin{equation*}
 \cD\bl((1+T^*T)^{1/2}\br)\supset\cD(T).
\end{equation*}

Now we see that, for a closed densely defined operator $T$, the operator $T(1+T^*T)^{-1/2}$ is
bounded and one has
\begin{equation*}
    \im T(1+T^*T)^{-1/2}=\im T.
\end{equation*}
Moreover,
\begin{equation*}
    \Ker T(1+T^*T)^{-1/2}=\Ker T,
\end{equation*}
since we have the following chain of equivalent relations:
\begin{align*}
    T(1+T^*T)^{-1/2}u=0&\Longleftrightarrow
    v=(1+T^*T)^{-1/2}u\in \Ker T\\
    &\Longleftrightarrow
    u=(1+T^*T)^{1/2}v=v\in \Ker T.
\end{align*}
(Here we have used the fact that $\Ker T=\Ker T^*T$ and hence, by the functional calculus,
$f(T^*T)v=f(0)v$ for any $v\in\Ker T$.)

Applying this reasoning to the operator $T=A_j$, we obtain all
assertions of the lemma except for the fact that the commutators
$[a,B_j]$ are compact operators for any $a\in\cA$. First, let us
prove the compactness of the commutators for $j>0$. Since the
domain $\cD(A_j)$ is invariant under $a$, we have
\begin{equation}\label{commcomp1}
    [a,B_j]=[a,A_j](1+\Delta_j)^{-1/2}+A_j[a,(1+\Delta_j)^{-1/2}]
\end{equation}
(both sides are defined everywhere  on $\cH_j$). The first term on
the right-hand side in~\eqref{commcomp1} is compact as the product
of the bounded operator $[a,A_j]$ by the compact operator
$(1+\Delta_j)^{-1/2}$. (The compactness of the second operator
follows, say, from~\cite[Chap.~V, Theorem~3.49]{Kat1} and the fact
that $(1+\Delta_j)^{-1}$ is compact by assumption.) It remains to
prove that the operator $A_j[a,(1+\Delta_j)^{-1/2}]$ is compact.
Let us temporarily omit the subscript $j$ for brevity. We have
(e.g., see Carey-Phillips~\cite{CaPh1})
\begin{equation}\label{useful}
    (1+\Delta)^{-1/2}=\frac1\pi\int_0^\infty
    \la^{-1/2}(1+\Delta+\la)^{-1}\,d\la\,,
\end{equation}
where the integral converges in operator norm. Hence we have
\begin{multline}\label{long1}
    A[a,(1+\Delta)^{-1/2}]=A\frac1\pi\int_0^\infty
    \la^{-1/2}(1+\Delta+\la)^{-1}[\Delta,a](1+\Delta+\la)^{-1}\,d\la\\
    =\frac1\pi\int_0^\infty
    \la^{-1/2}A(1+\Delta+\la)^{-1}[A^*,a]A(1+\Delta+\la)^{-1}\,d\la\\
    +\frac1\pi\int_0^\infty
    \la^{-1/2}A(1+\Delta+\la)^{-1}A^*[A,a](1+\Delta+\la)^{-1}\,d\la\\
    +\frac1\pi\int_0^\infty
    \la^{-1/2}A(1+\Delta+\la)^{-1}[\wt A^*,a]\wt A(1+\Delta+\la)^{-1}\,d\la\\
    +\frac1\pi\int_0^\infty
    \la^{-1/2}A(1+\Delta+\la)^{-1}\wt A^*[\wt A,a](1+\Delta+\la)^{-1}\,d\la\,.
\end{multline}
Here we have denoted $\wt A=A_{j-1}^*$. The closed operator $A$
can be passed through the integrals, since all integrals in
question converge in norm, as we shall see shortly. Obviously,
\footnote{Here the second and fourth operators are of course
considered on their (dense) domains, which suffices for our aims,
since the factors on the right of these operators in respective
integrands in~\eqref{long1} have ranges contained in these
domains.}
\begin{multline*}
    \norm{A(1+\Delta+\la)^{-1/2}},\norm{(1+\Delta+\la)^{-1/2}A^*},\\
    \norm{\wt A(1+\Delta+\la)^{-1/2}},\norm{(1+\Delta+\la)^{-1/2}\wt A^*}<1,
\end{multline*}
and so in each of the four integrals on the right-hand side in~\eqref{long1} the integrand is
compact and can be estimated in norm as $\const\cdot\la^{-3/2}$ at infinity and
$\const\cdot\la^{-1/2}$ at zero. Thus, the integrals converge absolutely and give compact
operators.

It remains to consider the case $j=0$. The preceding argument
would be of no use here, since the operator $\Delta_0$ is not
assumed to have a compact resolvent. However,
\begin{equation*}
    T(1+T^*T)^{-1/2}=\ov{(1+TT^*)^{-1/2}T}
\end{equation*}
for a closed densely defined operator $T$, where the bar on the right-hand side stands for the
closure. Using this and the direct sum expansion~\eqref{dirsum}, we readily obtain
\begin{equation*}
    A_0(1+\Delta_0)^{-1/2}=\ov{(1+\Delta_1)^{-1/2}A_0},
\end{equation*}
and then a computation similar to~\eqref{long1} applies. The proof of Lemma~\ref{lemmaa} is
complete.
\end{proof}
Now consider the rolled-up operator
\begin{equation}\label{dio}
    \wt D=B+B^*+P_0+(-1)^{n+1}P_n
\end{equation}
for the complex~\eqref{unbd-2}
\begin{lemma}\label{lemmab}
One has
\begin{equation}\label{piss-off}
    P_+(\wt D)=P_+(D).
\end{equation}
\end{lemma}
\begin{proof}[Proof of Lemma~\rom{\ref{lemmab}}]
It suffices to prove that
\begin{equation}\label{to-be-proved}
    P_+(A+A^*)=P_+(B+B^*),
\end{equation}
since the desired identity~\eqref{piss-off} can be obtained from~\eqref{to-be-proved} by adding
$P_0$ (and $P_n$ if $n$ is odd) on both sides. Next,
\begin{equation*}
    (A+A^*)^2=AA^*+A^*A=\bigoplus_{j=0}^n\Delta_j,
\end{equation*}
since $A^2=A^{*2}=0$. (Here each Laplacian acts in the respective direct summand
in~\eqref{summan}.) It follows, with regard to Lemma~\ref{lemmaa}, that
\begin{equation}\label{by-the-way}
    B+B^*=(A+A^*)\bl(1+(A+A^*)^2\br)^{-1/2}.
\end{equation}
It remains to note that, by the functional calculus of self-adjoint operators (e.g., see~\cite[\S
VI.5.2]{Kat1} or \cite[Chap.~9]{RiNa1}), it is always true that
\begin{equation*}
    P_+(T)=P_+\bl(T(1+T^2)^{-1/2}\br)
\end{equation*}
for a self-adjoint operator $T$, since the function
$\la\mapsto\la/(1+\la)^2$ vanishes at zero, is positive for
positive $\la$, and is negative for negative $\la$. This proves
the lemma.
\end{proof}
Now the preceding results actually prove the theorem, since we have $[P_0]=[P_+(\wt D)]$ by
Lemma~\ref{lemmaa} and Theorem~\ref{34}, and it remains to apply the result~\eqref{piss-off} of
Lemma~\ref{lemmab}. The proof of the theorem is complete.
\end{proof}

\begin{remark}
We note possible generalizations of the results obtained. First,
there is no need, in fact, to assume that the
sequence~\eqref{unbd-1} is exact in all intermediate terms; it
suffices to require that is has finite-dimensional cohomology in
these terms. Then the corresponding finite rank projections define
trivial quantizations and hence the assertion of the main theorem
remains valid. Second, if we replace the condition that the
resolvents of $\Delta_j$ be compact by the condition that they
belong to some von Neumann--Schatten classes, then in the end we
obtain, for some $p$, an unbounded $p$-summable Fredholm module
determined by the operator $D$.
\end{remark}

\appendix\section{The construction due to Epstein--Melrose--Mendoza}

In this appendix, we describe the constructions of~\cite{HHH,EpMe1} in slightly more detail,
restricting ourselves to the case in which the contact manifold $X$ considered there is the
cosphere bundle of a smooth closed compact manifold $M$.

Grauert~\cite{Gra1} showed that the cosphere bundle $X=S^*M$ can
always be represented as the boundary of a pseudoconvex Stein
manifold diffeomorphic to the ball bundle $B^*M$ (the Grauert
tube). Then the Szeg\"o--Calder{\'o}n projection $S$ can be defined
as the projection in $L^2(X)$ onto the subspace of boundary values
of holomorphic functions in the Grauert tube. The manifold $X=\pa
B^*M$ bears the natural Kohn--Rossi complex (the boundary
$\dbb$-complex; e.g., see~\cite{FoKo1})
\begin{equation}\label{dbb}
\begin{CD}
 0@>>>C^\infty(X)@>\dbb>>C^\infty(X,\La^{0,1})
 @>\dbb>>C^\infty(X,\La^{0,2})@>\dbb>>\dotsm\\&&&&@>\dbb>>
 C^\infty(X,\La^{0,n})@>>>0\,,
\end{CD}
\end{equation}
where $n$ is the dimension of the manifold $M$. Here $\La^{0,q}$ is the bundle of $(0,q)$-forms
associated with the natural complex structure on the contact hyperplane distribution.

The complex~\eqref{dbb} is acyclic in all dimensions except for
$0$ and $n$. Moreover, its cohomology space in the zero term
exactly coincides with the range of the Szeg\"o-Calder{\'o}n
projection $S$. Note that if we ``roll up'' the
complex~\eqref{dbb} in a standard way (i.e., pass to the operator
$\dbb+\dbb^*$ in the direct sum of all spaces of the
complex\footnote{The adjoint operator is taken with respect to
some Hermitian metric on $X$.}), then the resulting self-adjoint
operator is closely related to the Dirac operator on $X$. This
relationship can be conveniently represented in the block
$2\times2$ expansion of this operator corresponding to the
expansion of the space $\La^{0,*}$ into the sum of even and odd
components:
\begin{multline}\label{dirac}
    \dirac=\begin{pmatrix}
      -i\nabla_t & \dbb+\dbb^* \\
      \dbb+\dbb^* & i\nabla_t \\
    \end{pmatrix}
    :C^\infty(X,\La^{0,\operatorname{odd}})
\oplus C^\infty(X,\La^{0,\operatorname{even}})
    \\ \lra
    C^\infty(X,\La^{0,\operatorname{odd}})
\oplus C^\infty(X,\La^{0,\operatorname{even}})\,,
\end{multline}
where $\nabla_t$ is the covariant derivative along a direction transversal to the contact
distribution (more precisely, along the Reeb vector field).

This suggests that the Kohn--Rossi complex can be used to prove
the equivalence of the quantizations given by the
Szeg\"o--Calder{\'o}n projection $S$ and the positive spectral
projection of the Dirac operator. The construction carried out
in~\cite{HHH,EpMe1} consists of two steps.

1) Let
\begin{equation*}
 \wt S:C^\infty(X,\La^{0,n})\lra C^\infty(X,\La^{0,n})
\end{equation*}
be the orthogonal projection on the cokernel of the last operator $\dbb$ in the Kohn--Rossi
complex. The operator
\begin{equation}\label{chudo}
    G=S+(-1)^{n+1}\wt S+\dbb+\dbb^*
\end{equation}
(where the projections $S$ and $\wt S$ act in the first and the
last component, respectively) is self-adjoint and invertible. In
accordance with the general results proved in~Sec.~2, the
quantization with the help of the positive spectral projection of
$G$ is equivalent to the quantization with the help of the
Szeg\"o--Calder{\'o}n projection. (In~\cite{HHH}, only the
equalities for the indices of operator obtained from the same
symbol by these two quantizations were proved.)

2) Now it remains to deform the operator
\begin{equation*}
    G=\begin{pmatrix}
      S & \dbb+\dbb^*\\
      \dbb+\dbb^*& (-1)^{n+1}\wt S \\
    \end{pmatrix}
\end{equation*}
continuously in the class of Fredholm operators to the
Dirac operator $\dirac$. This procedure is described
in sufficient detail in the cited papers, and we only
briefly mention the main points. First, the operator
$G$ \textit{is not pseudodifferential} (owing to the
presence of the projections $S$ and $\wt S$, which are
not pseudodifferential operators). However, on contact
manifolds there is a natural ``anisotropic'' version
of pseudodifferential calculus, known as the
Heisenberg calculus. More precisely, in this calculus
the differentiations along contact directions have the
weight $1$, and the differentiation along the
transversal (Reeb) direction has the weight~$2$. Since
the contact distribution is nonintegrable, the product
in the algebra of principal symbols proves to be
noncommutative. The principal symbol over each point
$x\in X$ is a pair of equivariant operators in the
standard irreducible representation of the Heisenberg
group in $L^2(\RR^n)$. In particular, the projections
$S$ and $\wt S$ turn out to be Heisenberg
pseudodifferential operators whose symbols are the
projections on the vacuum vector of the harmonic
oscillator in one component of the pair and are equal
to zero in the other component. The paper~\cite{HHH}
contains also an ``extended'' version of the
Heisenberg calculus, which includes Heisenberg
pseudodifferential operators as well as usual
pseudodifferential operators. The homotopy of $G$ to
the Dirac operator $\dirac$ is carried out in the
framework of the extended calculus.

\providecommand{\bysame}{\leavevmode\hbox to3em{\hrulefill}\thinspace}



\end{document}